\RequirePackage{ifpdf}
\ifpdf 
\documentclass[pdftex]{sigma}
\else
\documentclass{sigma}
\fi

\newcommand{\la}{\label}

\def\exdeclare#1#2#3\par{\begin{#1}\label{#2}{\rm #3}\end{#1}\rm}
\def\stdeclare#1#2#3\par{\begin{#1}\label{#2}{#3}\end{#1}}
\def\stwrite#1#2{#1 \stn{#2}}  \def\stn{\ref}

 \def\th{\stwrite{Theorem}}

\def\nnorm#1{ \|\,{#1}\,\|}

\def\x1{{\xi }_{xx}}
\def\x2{{\xi }_{yy}}
\def\x3{{\xi }_{xy}}

\def\e1{{\eta }_{xx}}
\def\e2{{\eta }_{yy}}
\def\e3{{\eta }_{xy}}

\def\la{\lambda}

\def\Om{\Omega}

\def\l1{{\lambda}_1}

\def\nb{\nabla}
\def\kd{\partial}

\def\bb{\begin{equation}}
\def\ee{\end{equation}}
\def\ba{\begin{array}}
\def\ea{\end{array}}

\def\nb{\nabla}

\begin{document}

\allowdisplaybreaks

\renewcommand{\thefootnote}{$\star$}

\renewcommand{\PaperNumber}{055}

\FirstPageHeading

\ShortArticleName{Pohozhaev and Morawetz Identities in Elastostatics and Elastodynamics}

\ArticleName{Pohozhaev and Morawetz Identities \\ in Elastostatics and Elastodynamics\footnote{This paper is a
contribution to the Special Issue ``Symmetry, Separation, Super-integrability and Special Functions~(S$^4$)''. The
full collection is available at
\href{http://www.emis.de/journals/SIGMA/S4.html}{http://www.emis.de/journals/SIGMA/S4.html}}}

\Author{Yuri BOZHKOV~$^\dag$ and Peter J.~OLVER~$^\ddag$}

\AuthorNameForHeading{Y. Bozhkov and P.J.~Olver}

\Address{$^\dag$~Instituto de Matem\'atica, Estatistica e Computa\c c\~ao Cient\'\i fica - IMECC,\\
\hphantom{$^\dag$}~Universidade Estadual de Campinas - UNICAMP, Rua S\'ergio Buarque de Holanda, 651,\\
\hphantom{$^\dag$}~13083-859 - Campinas - SP, Brasil}

\EmailD{\href{mailto:bozhkov@ime.unicamp.br}{bozhkov@ime.unicamp.br}}
\URLaddressD{\url{http://www.ime.unicamp.br/~bozhkov/}}

\Address{$^\ddag$~School of Mathematics, University of Minnesota, Minneapolis, MN  55455, USA}
\EmailD{\href{mailto:olver@umn.edu}{olver@umn.edu}}
\URLaddressD{\url{http://www.math.umn.edu/~olver/}}

\ArticleDates{Received February 01, 2011, in f\/inal form June 02, 2011;  Published online June 08, 2011}

\Abstract{We construct identities of Pohozhaev type, in the context of elastostatics and elastodynamics, by using the Noetherian approach. As an application, a non-existence result for forced semi-linear isotropic and anisotropic elastic systems is established.}


\Classification{35J50; 35J47; 35L51}

\section{Introduction}

 Identities of Pohozhaev type have been widely used in the theory of partial dif\/ferential equations, in particular for establishing non-existence results for large classes of forced elliptic boundary value problems and eigenvalue problems,  \cite{p1,p2,ps}. The purpose of this note is to obtain and apply analogous identities in elastostatics and elastodynamics, which have not (to the authors' knowledge) been developed to date.  Our approach will be based on a fundamental identity f\/irst introduced by Noether in her seminal paper \cite{noe} that connected symmetries of variational problems to conservation laws of their Euler--Lagrange equations.

As noted in \cite{ol}, the identities originally due to Pohozhaev, \cite{p1,p2}, owe their existence to Noether's identity. For classical solutions of the linear equation $\Delta u+\lambda u=0$ such an identity was obtained by Rellich in \cite{rel1}. Further, in \cite{rel2}, Rellich established an integral identity for a function  belonging to certain function spaces, without any reference to dif\/ferential equations it may satisfy. The Rellich identity has been generalized by Mitidieri, \cite{em0,em}, for a pair of functions. General Rellich-type identities on Riemannian manifolds have been recently established in \cite{ye2,ye3} by use of Noether's identity applied to conformal Killing vector f\/ields.

In \cite{p1}, Pohozhaev established an integral identity for solutions of the Dirichlet problem for the semilinear Poisson equation $\Delta u+\lambda f(u)=0$ in a bounded domain with homogeneous Dirichlet boundary condition. Later, for solutions of general Dirichlet problems, he obtained in \cite{p2} what is now called the Pohozhaev identity. Such identities became very popular after the paper of Pucci and Serrin,  \cite{ps}, where, on p.~683, the relation with the general Noetherian theory is mentioned. See also the earlier paper by Knops and Stuart, \cite{ks}, and remarks in the second author's 1986 book \cite{ol}.  The Noetherian approach to Pohozhaev's identities was further developed and applied in \cite{ye,ye1,reich,vdv}.
Additional applications of Rellich--Pohozhaev estimates to nonlinear elliptic theory can be found in \cite{sr,sk}, while applications to nonlocal problems appear in \cite{ds}. With regard to geometric applications, \cite{be,dr,dru,scho} develop a systematic approach to Pohozhaev-type obstructions for partial dif\/ferential equations invariant under the action of a conformal group. For a relation between the Lie point symmetries of the nonlinear Poisson equation on a (pseudo-) Riemannian manifold and its isometry and conformal groups see \cite{yi}.

In dynamical problems, the conformal invariance of the wave and Klein--Gordon equations was used by Morawetz, \cite{mor}, to establish several very useful integral identities.  These were applied by her and Strauss, \cite{str1,str}, to the study of the decay, stability, and scattering of waves in nonlinear media.  In the f\/inal section, we will generalize Morawetz' conformal identity to some dynamical systems governing waves in elastic media.  Applications of our identity to decay and scattering of elastic waves will be treated elsewhere.

In elastostatics, the independent variables $x\in {\mathbb{R} }^n$, for $n\geq 2$, represent reference body coordinates, while the dependent variables
 $u=u(x)= (u^1(x), \ldots,u^n(x))$ represent the deformation of the point $x$.  The independent variable $x$ will belong to a bounded or unbounded domain $ \Om \subseteq {\mathbb{R} }^n $ that has suf\/f\/iciently (piecewise) smooth boundary $ \partial \Om$. We use $\nu $ to denote the outward unit normal on $\kd\Om $.  For elastodynamics, we append an additional independent variable, $t$, representing the time, and so $u = u(t,x)$.   The partial
 derivatives of a smooth (vector) function~$u(x)$ are denoted by subscripts:
 \[u^k_i:=\frac{\kd u^k}{\kd x_i}, \qquad u^k_t:=\frac{\kd u^k}{\kd t}, \qquad  u^k_{ij}:=\frac{{\kd }^2 u^k}{\kd x_i\kd x_j}, \qquad \mathrm{etc.}\]
  The $n \times n$ spatial Jacobian matrix $\nabla u = (u^k_i)$ is known as the \emph{deformation gradient}.

 We shall consistently use the Einstein summation convention over
repeated indices, which always run from $1$ to $n$.    We assume that all considered func\-tions, vector f\/ields,
 tensors, func\-tio\-nals, etc. are suf\/f\/iciently smooth in order that all
 the derivatives we write exist in the classical sense.  When we say that a function is ``arbitrary'', we mean that it is a suf\/f\/iciently smooth function of its arguments def\/ined on the domain~$\Om$. Extensions of our results to more general solutions will then proceed on a case by case basis.

 \section{Noether's identity}

A vector f\/ield
\[
\mathbf{v}={\xi }^i(x,u)\frac{\kd }{\kd x_i} + {\phi }^i (x,u) \frac{\kd }{\kd u^i}
\]
on the space of independent and dependent variables induces a f\/low that can be interpreted as a (local) one-parameter group of transformations.  The vector f\/ield is known as the \emph{infinitesimal generator} of the f\/low, \cite{ol}.  For example, the particular vector f\/ield
\[
\mathbf{v}=a x_i\frac{\kd }{\kd x_i} + b\, u^i\frac{\kd }{\kd u^i},
\]
where $a$, $b$ are constant, generates the group of scaling transformations
\[
(x,u) \ \longmapsto \ \big(\lambda ^a x, \lambda ^b u\big).
\]

The action of the group on functions $u=f(x)$ by transforming their graphs induces an action on their derivatives.  The corresponding inf\/initesimal generator of the prolonged group action has the form
 \begin{gather} \label{v1}
\mathrm{pr}^{(1)}\mathbf{v}  =   \xi ^i(x,u) \frac{\kd }{\kd x_i} + {\phi }^i (x,u) \frac{\kd }{\kd u^i} + {\phi }^i_j (x,u,\nabla u) \frac{\kd }{\kd u^i_j},
\end{gather}
where
 \begin{gather} \label{p1}
{\phi }^i_j (x,u,\nabla u) = D_j \phi^i - (D_j \xi^k) u^i_k = \frac{\kd \phi^i}{\kd x_j} + \frac{\kd \phi^i}{\kd u^k}u^k_j - \frac{\kd \xi^k}{\kd x_j} u^i_k  - \frac{\kd \xi^k}{\kd u^l}u^l_j u^i_k ,
\end{gather}
and $D_j = {\kd }/{\kd x_j} + u^k_j {\kd }/{\kd u^k}$ denotes the total derivative with respect to~$x_j$.
See \cite{ol} for a proof of this formula, along with its extension to higher order derivatives.

For a f\/irst order Lagrangian $L(x,u,\nabla u)$, \emph{Noether's identity} reads
 \begin{gather} \label{Noether}
\mathrm{pr}^{(1)}\mathbf{v} (L) +LD_i{\xi }^i = \mathrm E_{i }(L)({\phi}^{i } - u_j^{i}{\xi }^j) +D_i\left[ L{\xi }^i +
 \frac{\kd L}{\kd u_i^{j }}({\phi }^{j } - u_s^{j }{\xi }^s)\right],
 \end{gather}
where $\mathrm E_i$ is the Euler operator or variational derivative with respect to $u^i$, \cite{ol}.  Once stated, the verif\/ication of the identity is a straightforward computation. In the following sections, we will investigate how to use Noether's identity  in the framework of elasticity, and apply the corresponding integral identities to establish non-existence results. The proofs are sketched, while the full details are left to the interested reader as exercises.

\section{Elastostatics}

We recall that the equilibrium equations for a homogeneous
isotropic linearly elastic medium in the absence of body forces arise from the
variational principle with Lagrangian
\[ L_0 (x,u,\nabla u) = \frac{1}{2} \mu \nnorm{ \nabla u}^2 +\frac{1}{2}(\mu +\la ) ( \nabla \cdot u)^2 =
\frac{1}{2} \mu \sum_{i,j=1}^{n} \big(u^{i}_{j}\big)^2 +\frac{1}{2}(\mu +\la )\left( \sum_{i=1}^{n} u^{i}_{i}\right)^2, \]
where the parameters $\la $ and $\mu $ are the \emph{Lam\'e moduli}. The squared norm of the deformation gradient matrix $\nabla u$ refers to the sum of the squares of its entries, while $\nabla \cdot u$ denotes the divergence of the deformation. The corresponding Euler--Lagrange equations are known as \emph{Navier's equations}:
\[ \mu \Delta u +(\mu +\la )\nb (\nb \cdot u)=0,\]
where the Laplacian $ \Delta $ acts component-wise on $u$.
Henceforth, we assume that $\mu >0$ and $\mu +\la >0$, thereby ensuring strong ellipticity and positive def\/initeness of the underlying elasticity tensor, \cite{gur,ole1}.

In this paper, we shall study boundary value problems for elastic bodies that are subject to a nonlinear body-force potential $F(u)$. Thus, we modify the preceding Lagrangian
 \begin{gather*}
 L(x,u,\nabla u) = \frac{1}{2} \mu \nnorm{ \nabla u}^2 +\frac{1}{2}(\mu +\la ) ( \nabla \cdot u)^2 - F(u),
 \end{gather*}
 where we assume, without loss of generality, that $F(0) = 0$.
The associated equilibrium Euler--Lagrange equations are
\begin{gather}
\label{2} \mu \Delta u +(\mu +\la ) \nb (\nb \cdot u) + f(u)=0,
\end{gather}
where $f_i(u) = \kd F/\kd u^i$ are the components of  the gradient of the body-force potential with respect to the dependent variables $u$.

 More generally, we consider Lagrangians of the form:
  \begin{gather}\label{3}
 L=\frac{1}{2}  C^{k l}_{i j} e^{i}_{k} e^{j}_{l} -F(u),
 \end{gather}
where again $F(0) = 0$, and
 \begin{gather*}
e = \frac12 \big(\nabla u + \nabla u^T\big), \qquad \mathrm{with\ components} \quad e^i_k = \frac{1}{2} \big(u^i_k + u^k_i\big),
 \end{gather*}
is the \emph{strain tensor}. The quadratic components in the Lagrangian (\ref3) model the stored energy of a general anisotropic linearly elastic medium, while $F(u)$ represents a nonlinear body-force potential.
 The \emph{elastic moduli} $C^{k l}_{i j} $ are assumed to be constant, satisfying
 \begin{gather}\label{5}
 C^{k l}_{i j} = C^{ i l}_{k j}= C^{ k j}_{i l}= C^{ l k}_{j i} .
 \end{gather}
 Thus in planar elasticity there are 6 independent elastic moduli, while in
three dimensions 21 independent moduli are required in general, \cite{gur}. Additional
symmetry restrictions stemming from the constitutive properties of the
elastic material may place additional constraints on the moduli.
We may also assume
\begin{gather}
 \label{pos} C^{k l}_{i j} a^{i}_{k} a^{j}_{l} \geq 0
 \end{gather}
 for any matrix $A = (a^p_q)$. The less restrictive \emph{Legendre--Hadamard condition} is that
\begin{gather*}
C^{k l}_{i j} v^{i} v^{j} w_{k} w_{l} > 0
\end{gather*}
for any rank one matrix $A = v \otimes w$.   The Euler--Lagrange equations associated with (\ref3) read
 \begin{gather}\label{4}
 C^{k l}_{i j} u^{j}_{kl} +f_i(u)=0.
 \end{gather}

In general, the most basic Pohozhaev-type identity is based on the associated Noether identity for the inf\/initesimal generator of an adroitly chosen scaling transformation group, \cite{ol}.

\begin{theorem}\label{Poh1} Let $\Omega $ be a bounded domain in ${\mathbb{R}}^n$. Then the classical solutions of \eqref{4}~-- that is $u \in C^2(\Om)\cap C^1(\bar{\Om})$ -- subject to homogeneous Dirichlet boundary
conditions on $\kd \Om $ satisfy the following Pohozhaev-type identity:
 \begin{gather}\label{7}
\int_{\Om } \left[\frac{n-2}{2} u^k f_k(u)-nF(u) \right]dx = -\frac{1}{2}\int_{\kd\Om } C^{k l}_{i j} u^{i}_{k} u^{j}_{l}(x,\nu )ds,
\end{gather}
where $\nu $ is the outward unit normal to $\kd\Om $ and $(\cdot,\cdot)$ is the Euclidean scalar product in $ {\mathbb{R}}^n$.
\end{theorem}

\begin{proof} We consider the one-parameter group of dilations
\[(x,u) \ \longmapsto \ (\la x,{\la }^{(2-n)/2}u) \]
with inf\/initesimal generator
\[ \mathbf{v} = x_i \frac{\kd }{\kd x_i} +\frac{2-n}{2}u^{i }\frac{\kd }{\kd u^{i }} .\]
According to (\ref{v1}), (\ref{p1}), the f\/irst order prolongation of this vector f\/ield is
 \[\mathrm{pr}^{(1)}\mathbf{v} = x_i \frac{\kd }{\kd x_i} +\frac{2-n}{2}u^{i }\frac{\kd }{\kd u^{i }} -
\frac{n}{2} u_i^{j }\frac{\kd }{\kd u_i^{j }} .\]
Then one easily sees that
\begin{gather}\label{v2}
\mathrm{pr}^{(1)}\mathbf{v} (L) +LD_i{\xi }^i =\frac{n-2}{2} u^k f_k(u)-nF(u).
\end{gather}
The identity (\ref{7}) now follows from the divergence theorem using (\ref{Noether}), (\ref{3}), (\ref{v2}), our assumption $F(0) = 0$, and the homogeneous Dirichlet boundary conditions, taking into account that, on $\kd\Om $,
\begin{gather}\label{v3} u^j_s{\nu }_i= u^j_i{\nu }_s.\end{gather}
See \cite[p.~683]{ps} for more details on the last point.
\end{proof}

For the sake of completeness, we specialize the general elastic Pohozhaev identity to the isotropic case of the forced Navier equations (\ref{2}) in $\Omega $:
\begin{gather*}
\int_{\Om } \left[\frac{n-2}{2} u^k f_k(u)-nF(u) \right]dx = -\frac{1}{2}\int_{\kd\Om }\left[ \frac{1}{2}\,\mu \nnorm{ \nabla u}^2 +\frac{1}{2}(\mu +\la ) ( \nabla \cdot u)^2\right] (x,\nu )ds,
\end{gather*}
again  subject to homogeneous Dirichlet boundary conditions on $\kd \Om $.

As a corollary, we obtain the following non-existence result. Recall that the domain $\Om $ is \emph{star-shaped} with respect to the origin if  $(x,\nu )\geq 0$ for any $x\in\kd\Om $.

\begin{theorem}\label{nonxe}  Suppose that $\Om $ is a star-shaped domain. Let the function
 \[F=F(s)=F(s_1,\ldots,s_n)\in C^1({\mathbb{R}}^n)\] satisfy the conditions $F(0)=0$ and
\begin{gather}\label{vv} \frac{n-2}{2} s^k \frac{\kd F}{\kd s^k}-nF(s)\geq 0,\qquad  i=1,\ldots,n, \end{gather}
for any $s\in {\mathbb{R}}^n$.
We also suppose that the equality in \eqref{vv} holds if and only if $s=0$.
Then there is no non-trivial classical solution of the potential systems \eqref{2}, \eqref{4},  subject to homogeneous Dirichlet boundary conditions. \end{theorem}

\begin{proof} This theorem follows easily from the identity $(\ref{7})$, taking into account the positivity requirement (\ref{pos}) and star-shapedness condition. Indeed, any classical solution of (\ref{4}) subject to homogeneous Dirichlet boundary
conditions on $\kd \Om $ must satisfy the identity~(\ref{7}). For~(\ref{pos}) and the star-shapedness condition $(x,\nu )\geq 0$ for any $x\in\kd\Om $, it follows that the right-hand side of~(\ref{7}) is non-positive. On the other hand, by~(\ref{vv}) the left-hand side of (\ref{7}) is positive unless $u=0$ in $\Om$. Hence $u=0$.\end{proof}

\section{Elastodynamics}

In this section, we turn our attention to hyperbolic elastodynamic systems of potential type:
\begin{gather}\label{8}
-u^{i}_{tt} +C^{k l}_{i j} u^{j}_{kl} +f_i(u)=0 \end{gather}
in $\mathbb{R}\times \Om$ with homogeneous Dirichlet boundary
conditions on $\mathbb{R}\times \kd \Om$. The corresponding Lagrangian is given by
\begin{gather}\label{9}
L=\frac{1}{2} C^{k l}_{i j} e^{i}_{k} e^{j}_{l}  -\frac{1}{2} u^{i}_{t} u^{i}_{t} -F(u)
=\frac{1}{2} C^{k l}_{i j} u^{i}_{k} u^{j}_{l}  -\frac{1}{2} u^{i}_{t} u^{i}_{t} -F(u),
\end{gather}
where the second expression follows from the requirements (\ref{5}) on the elastic moduli.

\begin{theorem}\label{Pohdy} The classical solutions of the problem \eqref{8} satisfy the following identity
\begin{gather}
 \frac{d}{d t} \int_{\Om } \left[  t E(u) + u^{i}_{t}u^{i}_{k}x^k  +
\frac{n-1}{2} u^i u^{i}_{t} \right]
dx   =\int_{\Om } \left[ \frac{n-1}{2}u^k f_k(u)-(n+1)F(u)  \right]dx  \nonumber\\
 \qquad{} +\int_{\kd\Om } \left[ \frac{1}{2} \left(C^{k l}_{i j} u^{i}_{k} u^{j}_{l}  +\frac{1}{2} u^{i}_{t}u^{i}_{t}\right) (x,\nu ) +t\; C^{k l}_{i j} u^{i}_{k} u^{j}_{l}  \right]ds, \label{m1}
\end{gather}
where
\begin{gather*}
E(u)=\frac{1}{2} \big(C^{k l}_{i j} e^{i}_{k} e^{j}_{l}  + u^{i}_{t}u^{i}_{t}\big) - F(u)
=\frac{1}{2} \big(C^{k l}_{i j} u^{i}_{k} u^{j}_{l}  + u^{i}_{t}u^{i}_{t}\big) - F(u)
\end{gather*}
is the energy density.\end{theorem}

\begin{proof} We introduce a vector f\/ield $\mathbf{v}$ which is the inf\/initesimal generator of the dilation group
\[(t,x,u) \ \longmapsto \ \big(\la t,\la x,{\la }^{(1-n)/2} u\big). \]
The f\/irst order prolongation of $\mathbf{v} $ is given by
\[\mathrm{pr}^{(1)}\mathbf{v}= t\frac{\kd }{\kd t}+x_i \frac{\kd }{\kd x_i} +
 \frac{1-n}{2}u^i\frac{\kd }{\kd u^i}
 -\frac{n+1}{2}u_t^i\frac{\kd }{\kd u_t^i} -\frac{n+1}{2}u_j^i\frac{\kd }{\kd u_j^i}. \]
As a result,
 \begin{gather}\label{v4}
  \mathrm{pr}^{(1)}\mathbf{v}(L) +LD_i{\xi }^i =\frac{n-1}{2} u^k f_k(u)-(n+1)F(u),
\end{gather}
where the Lagrangian $L$ is given by (\ref{9}). Then, after some algebraic manipulations, the identi\-ty~(\ref{m1}) follows from the Noether identity~(\ref{Noether}) combined with (\ref{9}), (\ref{v3}), (\ref{v4}), the homogeneous Dirichlet boundary conditions, and, f\/inally, the divergence theorem. \end{proof}

Let $\Om \subset {\mathbb{R} }^n$ be a ball of radius $R$ centered at the origin. If we assume that $u(t,x)$ decays suf\/f\/iciently rapidly as  $R = |x| \rightarrow\infty $, then the following conformal identity holds for the nonlinear hyperbolic system (\ref{8}) in $\mathbb{R}\times {\mathbb{R} }^n$:

\begin{corollary}\label{nonxdy} The classical solutions of the problem \eqref{8} in  $\mathbb{R}\times {\mathbb{R} }^n$ that decay rapidly at large distances satisfy the identity
\[
\frac{d}{d t} \int_{{\mathbb{R} }^n} \left[ t E(u) + u^{i}_{t}u^{i}_{k}x^k  +
\frac{n-1}{2} u^i u^{i}_{t} \right]
dx =\int_{{\mathbb{R} }^n } \left[ \frac{n-1}{2}u^k f_k(u)-(n+1)F(u)  \right]dx.  \]
\end{corollary}

We observe that this result generalizes Morawetz's dilational identity for nonlinear wave equations, \cite{mor,str1,str}, to elastodynamical systems.

Finally, we consider a nonlinear hyperbolic system of
so-called Hamiltonian type, \cite{ye1},
\begin{gather}
  - u^{i}_{tt} + C^{k l}_{i j} u^{j}_{kl} +H_{v^i}  =  0,\nonumber \\
- v^{i}_{tt} + C^{k l}_{i j} v^{j}_{kl} +H_{u^i} =  0,  \label{11}
\end{gather}
in $\mathbb{R}\times \Om$ with homogeneous Dirichlet boundary
conditions on $\mathbb{R}\times \kd \Om $. (The independent va\-riab\-le~$x$ must belong to an even dimensional space ${\mathbb{R}}^{2m}$.) For such systems, we obtain a generalization of Morawetz's conformal identity~\cite{str}.

\begin{theorem}\label{Pohham} The classical solutions of the problem \eqref{11} satisfy the following identity
\begin{gather}
 \frac{d}{d t} \int_{\Om } \left[ t E(u,v) + \big(x^k u^{j}_{k}v^{j}_{t}  + x^k v^{j}_{k}u^{j}_{t}\big) +
\frac{n-1}{2} \big(a u^j v^{j}_{t}+bv^j u^{j}_t\big) \right]
dx \nonumber\\
\qquad
=\int_{\Om } \left[ \frac{n-1}{2}\big(a u^k H_{u^k} +b v^k H_{v^k}\big)  \right]dx
\nonumber\\ \qquad \quad {} +\int_{\kd\Om } \left[ \big(C^{k l}_{i j} u^{i}_{k} v^{j}_{l} + u^{i}_{t} v^{i}_{t}\big) (x,\nu ) +t  C^{k l}_{i j} \big(u^{i}_{t} v^{j}_{l} {\nu }_k + v^{j}_{t} u^{i}_{k} {\nu }_l \big) \right]ds,\label{m3}
\end{gather}
where the constants $a$ and $b$ are such that $a+b=2$ and
\begin{gather*}
E(u,v)=C^{k l}_{i j} u^{i}_{k} v^{j}_{l}  + u^{i}_{t} v^{i}_{t}-H(u,v).
\end{gather*}
\end{theorem}

\begin{proof} In order to prove \th{Pohham}, we use the same scheme as in the preceding \th{Pohdy}. Namely,
we consider a vector f\/ield $\mathbf{v}$ which is the inf\/initesimal generator of the dilation group
\[(t,x,u,v) \ \longmapsto \ \big(\la t,\la x,{\la }^{a(1-n)/2} u,{\la }^{b(1-n)/2} v\big), \]
 where the constants $a$ and $b$ satisfy $a+b=2$.
Applying the f\/irst order prolongation
\begin{gather*}
 \mathrm{pr}^{(1)}\mathbf{v}= t\frac{\kd }{\kd t}+x_i \frac{\kd }{\kd x_i} +
 \frac{a (1-n)}2u^i\frac{\kd }{\kd u^i}
 +\frac{b (1-n)}2v^i\frac{\kd }{\kd v^i}+\left(\frac{a (1-n)}2-1\right)u_t^i\frac{\kd }{\kd u_t^i} \\
 \phantom{\mathrm{pr}^{(1)}\mathbf{v}=}{}+\left(\frac{a (1-n)}2-1\right)u_j^i\frac{\kd }{\kd u_j^i}+\left(\frac{b (1-n)}2-1\right)v_t^i\frac{\kd }{\kd v_t^i} +\left(\frac{b (1-n)}2-1\right)v_j^i\frac{\kd }{\kd v_j^i}
 \end{gather*}
to the Lagrangian
 \begin{gather}\label{v5}
 L= \frac{1}{2} C^{k l}_{i j} u^{i}_{k} v^{j}_{l}- u_t^i v_t^i - H(u,v)
 \end{gather}
yields
 \begin{gather}\label{v6}
  {\rm pr}^{(1)}\mathbf{v} (L) +LD_i{\xi }^i = \frac{n-1}{2}\big(a u^k H_{u^k} +b v^k H_{v^k}\big),
\end{gather}
when $a+b=2$. Then, after some additional work, the identity (\ref{m3}) follows from (\ref{v5}), (\ref{Noether}), (\ref{v6}), (\ref{v3}), the homogeneous Dirichlet boundary conditions, and the divergence theorem. \end{proof}

\begin{corollary}\label{nonxham} Let $a$, $b$ and $E$ be as in Theorem~{\rm \ref{Pohham}}. Then, provided $u$ and $v$ decay sufficiently rapidly  at large distances,
\begin{gather*}
 \frac{d}{d t} \int_{{\mathbb{R}}^{2m} } \left[ t E(u,v) + \big(x^k u^{j}_{k}v^{j}_{t}  + x^k v^{j}_{k}u^{j}_{t}\big) +
\frac{n-1}{2} \big(a u^j v^{j}_{t}+bv^j u^{j}_t\big) \right]
dx \\
\qquad{} =\int_{{\mathbb{R}}^{2m} } \left[ \frac{n-1}{2}\big(a u^k H_{u^k} +b v^k H_{v^k} \big) \right]dx.
\end{gather*}
\end{corollary}

Applications of these identities to the stability and scattering of waves in elastic media  will be developed elsewhere.

\section{Further directions}

We emphasize that, in order to obtain Pohozhaev and Morawetz-type identities in elastostatics and elastodynamics by the Noetherian approach developed in \cite{ye,ye1}, we have focussed our attention on dilations, which are particular cases of conformal transformations. Further variational identities associated with other variational symmetries remain to be investigated. In particular, it would be interesting to analyze the variational identity for the semilinear Navier equations that corresponds to the f\/irst order generalized symmetry
\[\mathbf{v} = \bigl[\mu u^i_j +(2\mu +\la ){\delta }^i_j u^k_k\bigr]\frac{\kd }{\kd u^j} \]
found in \cite{ole2}. In fact, the variational and (at least in three dimensions) non-variational symmetries for isotropic linear elastostatics were completely classif\/ied in \cite{ole1,ole2} and the systems not only admit point symmetries, but also a number of f\/irst order generalized symmetries.  In the two-dimensional case, complex variable methods, as in \cite{musk}, are used to produce inf\/inite families of symmetries and conservation laws. Also in the two-dimensional case, additional symmetries appear when $3\mu + \lambda = 0$. In the three-dimensional case, when $7\mu + 3 \lambda =0$, Navier's equations admit a full conformal symmetry group, along with additional conformal-like generalized symmetries. Although these restrictions are non-physical, they still lead to interesting divergence identities in the more general isotropic case, which can be applied to the analysis of eigenvalue problems, and also, potentially, the nonlinearly forced case.  This remains to be investigated thoroughly.

\subsection*{Acknowledgements}

We wish to thank the referees for their useful suggestions. Yuri Bozhkov would also like to thank FAPESP and CNPq, Brasil, for partial f\/inancial support. Peter Olver was supported in part by NSF Grant DMS 08--07317.  We both would like to thank FAPESP, S\~ao Paulo, Brasil, for the grant giving Peter Olver the opportunity to visit IMECC-UNICAMP, where this work was initiated.

\pdfbookmark[1]{References}{ref}
\LastPageEnding


\begin{thebibliography}{99}

\footnotesize\itemsep=0pt

\bibitem{be}
Bourguignon J.-P.,   Ezin J.-P.,
Scalar curvature functions in a conformal class of metrics and conformal transformations,
\href{http://dx.doi.org/10.2307/2000667}{{\it Trans. Amer. Math. Soc.}} {\bf 301} (1987),   723--736.

\bibitem{yi}
Bozhkov  Y.,  Freire  I.L.,
Special conformal groups of a Riemannian manifold and Lie point symmetries of the nonlinear Poisson equation,
\href{http://dx.doi.org/10.1016/j.jde.2010.04.011}{{\it J.~Differential Equations}} {\bf 249} (2010), 872--913,
\href{http://arxiv.org/abs/0911.5292}{arXiv:0911.5292}.

\bibitem{ye}
Bozhkov Y.,   Mitidieri E.,
The Noether approach to Pohozhaev's identities,
\href{http://dx.doi.org/10.1007/s00009-007-0125-y}{{\it Mediterr. J. Math.}} {\bf 4} (2007), 383--405.

\bibitem{ye1}
Bozhkov Y.,   Mitidieri E.,
Lie symmetries and criticality of semilinear  dif\/ferential systems,
\href{http://dx.doi.org/10.3842/SIGMA.2007.053}{{\it SIGMA}} {\bf 3} (2007),   053, 17~pages,
\href{http://arxiv.org/abs/math-ph/0703071}{math-ph/0703071}.

\bibitem{ye2}
Bozhkov Y.,  Mitidieri  E.,
Conformal Killing vector f\/ields and  Rellich type identities on Riemannian manifolds.~I,
in Geometric Methods in PDE's, {\it Lect. Notes Semin. Interdiscip. Mat.}, Vol.~7, Semin. Interdiscip. Mat. (S.I.M.), Potenza, 2008, 65--80.


\bibitem{ye3}
Bozhkov Y.,   Mitidieri  E.,
 Conformal Killing vector f\/ields and  Rellich type identities on Riemannian manifolds.~II,
 \href{http://dx.doi.org/10.1007/s00009-011-0126-8}{{\it Mediterr. J. Math.}}, to appear,
\href{http://arxiv.org/abs/1012.2993}{arXiv:1012.2993}.

 \bibitem{dr}
 Delano\"e  P.,  Robert  F.,
 On the local Nirenberg problem for the $Q$-curvatures,
 \href{http://dx.doi.org/10.2140/pjm.2007.231.293}{{\it Pacific~J. Math.}} {\bf 231} (2007), 293--304,
\href{http://arxiv.org/abs/math.DG/0601732}{math.DG/0601732}.

\bibitem{ds}
Dolbeault  J.,  Sta{\'n}czy  R.,
Non-existence and uniqueness results for supercritical semilinear elliptic equationss,
\href{http://dx.doi.org/10.1007/s00023-009-0016-9}{{\it Ann. Henri Poincar\'e}} {\bf 10} (2010), 1311--1333,
\href{http://arxiv.org/abs/0901.0224}{arXiv:0901.0224}.

 \bibitem{dru}
 Druet O.,
 From one bubble to several bubbles: the low-dimensional case,
{\it J. Differential Geom.} {\bf 63} (2003), 399--473.

 \bibitem{gur}
 Gurtin M.E.,
 The linear theory of elasticity, in Handbuch der Physik, Vol. VIa/2, Editor C.~Truesdell, Springer-Verlag, New York, 1972,   1--295.

 \bibitem{ks}
 Knops R.J.,  Stuart C.A.,
 Quasiconvexity and uniqueness of equilibrium solutions in nonlinear elasticity,
 \href{http://dx.doi.org/10.1007/BF00281557}{{\it Arch. Rational Mech. Anal.}} {\bf 86} (1984), 233--249.

\bibitem{em0}
Mitidieri E.,
 A Rellich identity and  applications,
Rapporti Interni No~25, Univ. Udine, 1990, 35~pages.

\bibitem{em}
Mitidieri  E.,
A Rellich type identity and  applications,
\href{http://dx.doi.org/10.1080/03605309308820923}{{\it Comm. Partial Differential Equations}} {\bf  18} (1993), 125--151.

\bibitem{mor}
Morawetz C.S.,
Notes on time decay and scattering for some hyperbolic problems,
{\it Regional Conference Series in Applied Mathematics}, no.~19, Society for Industrial and Applied Mathematics, Philadelphia, Pa., 1975.

\bibitem{musk}
Muskhelishvili N.I.,
Some basic problems of the mathematical theory of elasticity. Fundamental equations, plane theory of elasticity, torsion and bending, Noordhof\/f, Groningen, 1953.

\bibitem{noe}
Noether E., Invariante Variationsprobleme,
{\it Nachr. Konig. Gesell. Wissen. G\"ottingen, Math.-Phys. Kl.} (1918), 235--257.

\bibitem{ole1}
Olver P.J.,
Conservation laws in elasticity. I.~General results,
\href{http://dx.doi.org/10.1007/BF00281447}{{\it Arch. Rat. Mech. Anal.}} {\bf 85} (1984), 111--129.

\bibitem{ole2}
Olver P.J., Conservation laws in elasticity. II.~Linear homogeneous isotropic elastostatics,
\href{http://dx.doi.org/10.1007/BF00281448}{{\it Arch. Rat. Mech. Anal.}} {\bf 85} (1984), 131--160,
 Errata, \href{http://dx.doi.org/10.1007/BF00251537}{{\it Arch. Rat. Mech. Anal.}} {\bf 102} (1988), 385--387.

 \bibitem{ol}
 Olver P.J.,
Applications of Lie groups to dif\/ferential equations, 2nd ed., {\it Graduate Texts in Mathematics}, Vol.~107,
Springer-Verlag, New York, 1993.

 \bibitem{p1}
 Pohozhaev S.I.,
 On the eigenfunctions of the equation $\Delta u+\lambda f(u)=0$,
 {\it Dokl. Akad. Nauk SSSR} {\bf 165} (1965), 36--39
 (English transl.: {\it Soviet Math. Dokl.} {\bf 6} (1965), 1408--1411).

 \bibitem{p2}
 Pohozhaev S.I.,
 On eigenfunctions of quasilinear elliptic problems,
 {\it Mat. Sb.} {\bf 82} (1970), 192--212  (English transl.: \href{http://dx.doi.org/10.1070/SM1970v011n02ABEH001133}{{\it Math. USSR Sbornik}} {\bf 11} (1970), 171--188).

 \bibitem{ps}
 Pucci P., Serrin J.,
 A general variational identity,
 \href{http://dx.doi.org/10.1512/iumj.1986.35.35036}{{\it Indiana Univ. Math.~J.}} {\bf 35} (1986), 681--703.

 \bibitem{reich}
 Reichel W.,
 Uniqueness theorems for variational problems by the method of transformation groups,
 {\it Lecture Notes in Mathematics}, Vol.~1841, Springer-Verlag, Berlin, 2004.

 \bibitem{rel1}
 Rellich F.,
 Darstellung der Eigenverte von $\Delta u+\lambda u=0$ durch ein Randintegral,
 \href{http://dx.doi.org/10.1007/BF01181459}{{\it Math.~Z.}} {\bf 46} (1940), 635--636.

 \bibitem{rel2}
 Rellich F.,
 Halbbeschr\"ankte Dif\/ferentialoperatoren h\"oherer Ordnung,
 in  Proceedings of the International Congress of Mathematicians (Amsterdam, 1954),
Vol.~III, Erven P.~Noordhof\/f N.V., Groningen; North-Holland Publishing Co., 1956, 243--250.

 \bibitem{sr}
 Schaaf R.,
 Uniqueness for semilinear elliptic problems: supercritical growth and domain geometry,
 {\it Adv. Differential Equations} {\bf 5} (2000), 1201--1220.

\bibitem{sk}
Schmitt K.,
Positive solutions of semilinear elliptic boundary value problems,
in Topological Methods in Dif\/ferential Equations and Inclusions (Montreal, PQ, 1994),
{\it NATO Adv. Sci. Inst. Ser. C Math. Phys. Sci.}, Vol.~472, Kluwer Acad. Publ., Dordrecht, 1995,  447--500.

\bibitem{scho}
Schoen R.M.,
Variational theory for the total scalar curvature functional for Riemannian metrics and related topics,
 in Topics in Calculus of Variations (Montecatini Terme, 1987), {\it Lecture Notes in Math.}, Vol.~1365, Springer, Berlin, 1989, 120--154.

\bibitem{str1}
Strauss W.A.,
Nonlinear invariant wave equations,
in Invariant Wave Equations, Editors G.~Velo and A.S.~Wightman,
{\it Lecture Notes in Physics}, Vol.~73, Springer-Verlag, New York, 1978,  197--249.

\bibitem{str}
Strauss W.A.,
Nonlinear wave equations, {\it CBMS Regional Conference Series}, Vol.~73, Amer. Math. Soc., Providence, R.I., 1989.

\bibitem{vdv}
van der Vorst R.C.A.M.,
Variational identities and applications to dif\/ferential systems,
\href{http://dx.doi.org/10.1007/BF00375674}{{\it Arch. Rat. Mech. Anal.}} {\bf 116} (1991), 375--398.

\end{thebibliography}
\end{document}